\documentclass[12pt]{article}
\usepackage{latexsym}
\usepackage{amssymb}
\usepackage{graphicx}

\newtheorem{Theorem}{Theorem}[part]
\newtheorem{Definition}{Definition}[part]
\newtheorem{Proposition}{Proposition}[part]

\newtheorem{Lemma}{Lemma}[part]
\newtheorem{Corollary}{Corollary}[part]

\newtheorem{Example}{Example}[part]

\def \ep{\hbox{ }\hfill$\Box$}

\begin{document}

\title{Global Uniqueness and Solvability for Tensor Complementarity Problems \thanks{This work is partially supported by the National Natural Science Foundation of
China (Grant No. 11171252 and Grant No. 11431002).}}

\author{Xue-Li Bai\thanks{Department of Mathematics, School of Science, Tianjin University, Tianjin 300072, P.R. China. Email: bai7707@163.com}
\and
Zheng-Hai Huang\thanks{Corresponding Author. Department of Mathematics, School of Science, Tianjin University, Tianjin 300072, P.R. China. This author is also with the Center for Applied Mathematics of Tianjin University. Email: huangzhenghai@tju.edu.cn. Tel:+86-22-27403615 Fax:+86-22-27403615}
\and
Yong Wang\thanks{Department of Mathematics, School of Science, Tianjin University, Tianjin 300072, P.R. China. Email: wang\_yong@tju.edu.cn. }
}

\date{July 26, 2015}

\maketitle

\begin{abstract}
Recently, the tensor complementarity problem (TCP for short) has been investigated in the literature. An important question involving the property of global uniqueness and solvability (GUS-property) for a class of TCPs was proposed by Song and Qi in their paper ``Properties of Some Classes of Structured Tensors". In the present paper, we give an answer to this question by constructing two counter-examples. We also show that the solution set of this class of TCPs is nonempty and compact. In particular, we introduce a class of related structured tensors, and show that the corresponding TCP has the GUS-property. \vspace{3mm}

\noindent {\bf Key words:}\hspace{2mm} Tensor complementarity problem; Nonlinear complementarity problem; Global uniqueness and solvability; $P$ tensor \and Strong $P$ tensor\vspace{3mm}

\noindent {\bf Mathematics Subject Classifications(2000):}\hspace{2mm} 90C33; 65K10; 15A18; 15A69; 65F15; 65F10
\end{abstract}

\section{Introduction}
\label{intro}

As a generalization of the linear complementarity problem \cite{cps92}, the tensor complementarity problem has been introduced and investigated in \cite{sq14a,sq14b,sq15,cqw15,lqx15}, which is a specific class of nonlinear complementarity problems \cite{hp90,fp03,hxq06}. By using properties of structured tensors, many good results for the tensor complementarity problem have been obtained in the literature.

It is well-known that a linear complementarity problem has the property of global uniqueness and solvability (GUS-property) if and only if the matrix involved in the concerned problem is a $P$-matrix \cite{stw58}. It is also well-known that such a result can not be generalized to the nonlinear complementarity problem \cite{fp03,hxq06}. {\it A natural question is whether such a result can be generalized to the tensor complementarity problem or not}? i.e., whether the result {\it that a tensor complementarity problem has the GUS-property if and only if the tensor involved in the problem is a $P$ tensor} holds or not? Such a question was proposed by Song and Qi (see, Question 6.3 in \cite{sq14a}). In this paper, we show that the answer to this question is negative by constructing two counter-examples.

It has been shown that the tensor complementarity problem with a $P$ tensor has a solution; while our counter-example demonstrates that it is possible such a problem has more than one solution. Thus, a natural question is what more we can say about the solution set of such a complementarity problem; and another natural question is that for which kind of tensor, the corresponding tensor complementarity problem has the GUS-property. For the first question, we will show that the solution set of the tensor complementarity problem is nonempty and compact when the involved tensor is a $P$ tensor; and for the second question, we will introduce a new class of tensors, called  the {\it strong $P$ tensor}, and show that the corresponding tensor complementarity problem has the GUS-property. We also show that the set consists of all strong $P$ tensors is a proper subset of the set  consists of all $P$ tensors; and hence, many results obtained for the case of the $P$ tensor are still satisfied for the case of the strong $P$ tensor.

The rest of this paper is organized as follows. In the next section, we first briefly review some basic concepts and results which are useful in the subsequent analysis. In Sect. 3, we give a negative answer to the result that a tensor complementarity problem has the GUS-property if and only if the tensor involved in the problem is a $P$ tensor, and show that the solution set of the tensor complementarity problem with a $P$ tensor is compact. In Sect. 4, we introduce the concept of the strong $P$ tensor and discuss its related properties. Conclusions are given in Sect. 5.

\section{Preliminaries}
\label{sec:1}

The complementarity problem, denoted by CP$(F)$, is to find a point $x \in \mathbb{R}^n$ such that
\begin{eqnarray*}
x\geq 0,\quad F(x)\geq 0,\quad x^TF(x)=0.
\end{eqnarray*}
When $F(x)=Ax+q$ with given $A\in \mathbb{R}^{n\times n}$ and $q\in \mathbb{R}^n$, CP$(F)$ reduces to the linear complementarity problem (denoted by LCP$(q,A)$), which is to find a point $x \in \mathbb{R}^n$ such that
\begin{eqnarray*}
x\geq 0,\quad Ax+q\geq 0,\quad x^T(Ax+q)=0;
\end{eqnarray*}
and when $F(x)={\cal A}x^{m-1}+q$ with given ${\cal A}=(a_{i_1i_2\cdots i_m})\in T_{m,n}$ (the set of all real $m$th order $n$-dimensional tensors) and $q\in \mathbb{R}^n$, CP$(F)$ reduces to the tensor complementarity problem (denoted by TCP$(q,{\cal A})$) \cite{sq14a,sq14b,sq15,cqw15,lqx15}, which is to find a point $x \in \mathbb{R}^n$ such that
\begin{eqnarray*}
x\geq 0,\quad {\cal A}x^{m-1}+q\geq 0,\quad x^T({\cal A}x^{m-1}+q)=0,
\end{eqnarray*}
where ${\cal A}x^{m-1}\in \mathbb{R}^n$ defined by
\begin{eqnarray*}
({\cal A}x^{m-1})_i:=\sum_{i_2,\cdots,i_m=1}^na_{ii_2\cdots i_m}x_{i_2}\cdots x_{i_m},\quad \forall i\in \{1,2,\ldots,n\}.
\end{eqnarray*}
It is easy to see that
\begin{eqnarray*}
{\cal A}x^m=x^T({\cal A}x^{m-1})=\sum_{i_1, i_2, \cdots , i_m = 1}^na_{i_1i_2\cdots i_m}x_{i_1}x_{i_2}\cdots x_{i_m}.
\end{eqnarray*}

Throughout this paper, for any positive integer $n$, we denote $[n]:=\{1,2,\ldots,n\}$ and $\mathbb{R}_+^n:=\{x\in \mathbb{R}^n: x\geq 0\}$. For any $x\in \mathbb{R}^n$, we denote
$$
[x]_+:=(\max\{x_1,0\},\ldots,\max\{x_n,0\})^T.
$$

The eigenvalue of tensor is initially studied by Qi \cite{q05} and Lim \cite{l05}. If there is a nonzero vector $x\in \mathbb{R}^n$ and a scalar $\lambda\in \mathbb{R}$ such that
$$
({\cal A}x^{m-1})_i=\lambda x_i^{m-1},\quad \forall i\in [n],
$$
then $\lambda$ is called an $H$-eigenvalue of ${\cal A}$ and $x$ is called an $H$-eigenvector of ${\cal A}$ associated with $\lambda$; and if there is a nonzero vector $x\in \mathbb{R}^n$ and a scalar $\lambda\in \mathbb{R}$ such that
$$
{\cal A}x^{m-1}=\lambda x,\;\forall i\in \{1,2,\ldots,n\}\quad\mbox{\rm and}\quad x^Tx=1,
$$
then $\lambda$ is called a $Z$-eigenvalue of ${\cal A}$ and $x$ is called a $Z$-eigenvector of ${\cal A}$ associated with $\lambda$.

Recently, many classes of structured tensors are introduced and the related properties are studied \cite{sq14a,sq14b,sq15,cqw15,lqx15,q13,zqz14,cq15,qs14,dqw13,sq14c,cq14}. In this paper, we need the following concepts of several structured tensors.

\begin{Definition}\label{def-PT}
Let ${\cal A}=(a_{i_1\cdots i_m})\in T_{m,n}$. We say that ${\cal A}$ is
\begin{itemize}
\item[(i)] a {\bf strictly semi-positive tensor} iff for each $x\in \mathbb{R}^n_+\setminus \{0\}$, there exists an index $i\in [n]$ such that
    $x_i>0$ and $({\cal A}x^{m-1})_i>0$;
\item[(ii)]  a {\bf $P$ tensor} iff for each $x\in \mathbb{R}^n\setminus \{0\}$, there exists an index $i\in [n]$ such that $x_i({\cal A}x^{m-1})_i>0$;
\item[(iii)] an {\bf $R$-tensor} iff there is no $(x,t)\in (\mathbb{R}^n_+\setminus\{0\})\times \mathbb{R}_+$ such that for any $i\in [n]$,
\begin{eqnarray*}
\left\{\begin{array}{l}
({\cal A}x^{m-1})_i+t=0\; \mbox{\rm if}\; x_i>0,\\
({\cal A}x^{m-1})_i+t\geq 0\; \mbox{\rm if}\; x_i=0.
\end{array}\right.
\end{eqnarray*}
\end{itemize}
\end{Definition}

Obviously, every $P$ tensor is a strictly semi-positive tensor and an $R$-tensor.

In this paper, we also need the following concepts of functions.

\begin{Definition}[\cite{sq14b}]\label{def-PF}
Let mapping $F: K\subseteq \mathbb{R}^n \rightarrow \mathbb{R}^n$. We say that $F$ is
\begin{itemize}
\item[(i)] a {\bf $P$-function} if for all pairs of distinct vectors $x$ and $y$ in $K$,
\begin{eqnarray*}
\max_{i \in [n]}\;(x_i-y_i)(F_i(x)-F_i(y)) > 0;
\end{eqnarray*}
\item[(ii)] a {\bf uniform $P$-function} if there exists a constant $\mu>0$ such that for all pairs of vectors $x$ and $y$ in $K$,
\begin{eqnarray*}
\max_{i \in [n]}\;(x_i-y_i)(F_i(x)-F_i(y)) \geq \mu\|x-y\|^2.
\end{eqnarray*}
\end{itemize}
\end{Definition}

Obviously, every uniform $P$-function is a $P$-function. In addition, it is easy to see from Definition \ref{def-PT} and Definition \ref{def-PF} that if the mapping ${\cal A}x^{m-1}+q$ with any given $q\in \mathbb{R}^n$ is a $P$-function, then ${\cal A}$ is a $P$ tensor.

Recall that LCP$(q,A)$ is said to have {\bf the GUS-property} if LCP$(q,A)$ has a unique solution for every $q\in \mathbb{R}^n$. Similarly, we say that {\it TCP$(q,{\cal A})$ has the GUS-property if TCP$(q,{\cal A})$ has a unique solution for every $q\in \mathbb{R}^n$}. For the solvability of LCP$(q,A)$, an important result is that LCP$(q,A)$ has the GUS-property iff the matrix $A$ is a $P$-matrix. In fact, the GUS-property has been extensively discussed for various complementarity problems, including nonlinear complementarity problems \cite{mk77}, linear complementarity problems over symmetric cones \cite{gs07} and Lorentz cone linear complementarity problems on Hilbert spaces \cite{mh11}. A natural question is given by
\begin{itemize}
\item[{\bf Q1}] {\it Whether or not TCP$(q,{\cal A})$ has the GUS-property iff the tensor ${\cal A}$ is a $P$ tensor}?
\end{itemize}
Such a question was proposed by Song and Qi (see, Question 6.3 in \cite{sq14a}). In the next section, we will answer this question.

\section{Answer to {\bf Q1}}
\label{sec:2}

First, we construct a TCP$(q,{\cal A})$ which has a unique solution for every $q\in \mathbb{R}^2$.
\begin{Example}\label{exam-1}
Let ${\cal A}=(a_{i_1 i_2 i_3})\in T_{3,2}$, where $a_{111}=1,a_{222}=1$ and all other $a_{i_{1}i_{2}i_{3}}=0$. Then,
$$
{\cal A}x^2=\left(\begin{array}{c}
x_{1}^2 \\ x_{2}^2
\end{array}\right).
$$
In this case, TCP$(q,{\cal A})$ is to find $x\in \mathbb{R}^2$ such that
\begin{eqnarray}\label{E-exam-1-1}
\left\{\begin{array}{l}
x_1 \geq 0,\\
x_2 \geq 0,
\end{array}\right.\quad
\left\{\begin{array}{l}
x_{1}^2+q_1 \geq 0,\\
x_{2}^2+q_2 \geq 0,
\end{array}\right. \quad\mbox{\rm and}\quad
\left\{\begin{array}{l}
x_1(x_{1}^2+q_1) = 0,\\
x_2(x_{2}^2+q_2) = 0.
\end{array}\right.
\end{eqnarray}
For any $q\in \mathbb{R}^2$, let $x^q:=(x_1^q,x_2^q)\in \mathbb{R}^2$ be given by
$$
x_i^q:=\left\{\begin{array}{ll}
0,\quad & \mbox{if}\; q_i\geq 0,\\
\sqrt{-q_i},\quad & \mbox{\rm otherwise}.
\end{array}\right.
$$
It is easy to see that for every $q\in \mathbb{R}^2$, TCP$(q,{\cal A})$ given by (\ref{E-exam-1-1}) has a unique solution $x^q$. Thus, we obtain that TCP$(q,{\cal A})$ given in this example has the GUS-property.
\end{Example}

It is proved in \cite[Proposition 2.1]{yy14} that there dose not exist an odd order $P$ tensor; and hence, the tensor given in Example \ref{exam-1} is not a $P$ tensor. This, together with Example \ref{exam-1}, implies that one does not obtain that tensor ${\cal A}$ is a $P$ tensor under the assumption that TCP$(q,{\cal A})$ has the GUS-property.

Second, we construct the following TCP$(q,{\cal A})$ where ${\cal A}\in T_{4,2}$ is a $P$ tensor, but it has two distinct solutions for some $q\in \mathbb{R}^2$.
\begin{Example}\label{exam-2}
Let ${\cal A}=(a_{i_1 i_2 i_3 i_4})\in T_{4,2}$, where $a_{1111}=1,a_{1112}=-2,a_{1122}=1,a_{2222}=1$ and all other $a_{i_{1}i_{2}i_{3}i_{4}}=0$. Then,
$$
{\cal A}x^3=\left(\begin{array}{c}
x_{1}^3-2x_1^2x_2+x_1x_2^2 \\ x_{2}^3
\end{array}\right),
$$
and
$$
x_1({\cal A}x^{3})_1=x_1^4-2x_1^3x_2+x_1^2x_2^2,\quad x_2({\cal A}x^{3})_2=x_2^4.
$$
For any $x\in \mathbb{R}^2\setminus\{0\}$, it is easy see that
\begin{itemize}
\item when $x_2 \neq 0$, it follows that $x_2({\cal A}x^{3})_2 > 0$; and
\item when $x_2 = 0$, it follows that $x_1 \neq 0$ since $x \neq 0$, and in this case, we have $x_1({\cal A}x^{3})_1=x_1^4 >0$.
\end{itemize}
Thus, for any $x \in \mathbb{R}^2\setminus\{0\}$, there is at least one index $i \in \{1,2\}$ such that $x_i({\cal A}x^{3})_i>0$. So, we obtain that tensor ${\cal A}$ given in this example is a $P$ tensor by Definition \ref{def-PT}(ii).

Taking $q=(0,-1)^T$, we consider TCP$(q,{\cal A})$ of finding $x\in \mathbb{R}^2$ such that
\begin{eqnarray}\label{E-exam-2-1}
\left\{\begin{array}{l}
x_1 \geq 0,\\
x_2 \geq 0,
\end{array}\right.\;\;
\left\{\begin{array}{l}
x_{1}^3-2x_1^2x_2+x_1x_2^2 \geq 0,\\
x_{2}^3-1 \geq 0,
\end{array}\right.\;\; \mbox{\rm and}\;\;
\left\{\begin{array}{l}
x_1(x_{1}^3-2x_1^2x_2+x_1x_2^2) = 0,\\
x_2(x_{2}^3-1) = 0.
\end{array}\right.
\end{eqnarray}
It is easy to see that both $x=(0,1)^T$ and $x=(1,1)^T$ are the solutions to TCP$(q,{\cal A})$ given by (\ref{E-exam-2-1}).
\end{Example}

From Example \ref{exam-2}, we obtain that TCP$(q,{\cal A})$ with ${\cal A}$ being a $P$ tensor does not possess the GUS-property.

Note that TCP$(q,{\cal A})$ with ${\cal A}$ being a $P$ tensor has a solution for every $q\in \mathbb{R}^n$ by \cite[Corollary 3.3]{sq14b}, and that it is possible that this class of complementarity problems has more than one solution from Example \ref{exam-2}. What more can we say about the solution set of this class of complementarity problems? In the following, we show that the solution set of TCP$(q,{\cal A})$ with ${\cal A}$ being a $P$ tensor is compact.
\begin{Theorem}\label{the-ncss}
For any $q \in \mathbb{R}^n$ and a $P$ tensor ${\cal A}\in T_{m,n}$, the solution set of TCP$(q,{\cal A})$ is nonempty and compact.
\end{Theorem}

{\bf Proof.}
Since ${\cal A}$ is a $P$ tensor, it follows from \cite[Corollary 3.3]{sq14b} that TCP$(q,{\cal A})$ has a solution for every $q\in \mathbb{R}^n$. So we only need to show that the solution set of TCP$(q,{\cal A})$ is compact. We divide the proof into two parts.

{\bf Part 1}. We show the boundedness of the solution set. To this end, we first show the following result:
\begin{itemize}
\item[{\bf R1}] If there is a sequence $\{ x^k \} \subset \mathbb{R}_+^n$ satisfying
\begin{eqnarray}\label{c0}
\| x^k \| \rightarrow \infty\quad \mbox{\rm and}\quad \frac{[-{\cal A}(x^k)^{m-1}-q]_+}{\| x^k \|} \rightarrow 0\;\; \mbox{\rm as}\; k \rightarrow \infty,
\end{eqnarray}
then there exists an $i \in [n]$ such that $x_{i}^{k}[{\cal A}(x^k)^{m-1}+q]_i > 0$ holds for some $k \geq 0$.
\end{itemize}
In the following, we assume that the result {\bf R1} does not hold and derive a contradiction. Given an arbitrary sequence $\{ x^k \} \subset \mathbb{R}_+^n$ satisfying (\ref{c0}), then since the result {\bf R1} does not hold, we have that
\begin{eqnarray}\label{c0-1}
x_{i}^{k}[{\cal A}(x^k)^{m-1}+q]_i \leq 0\quad \forall i \in [n],\;\forall k \geq 0.
\end{eqnarray}
Since the sequence $\{ \frac{x^k}{\| x^k \|}\}$ is bounded, without loss of generality, we can assume $\lim_{k \rightarrow \infty}\frac{x^k}{\| x^k \|} = \bar{x} \in \mathbb{R}^n$. From $\{ x^k \} \subset \mathbb{R}_+^n$ and $\| x^k \| \rightarrow \infty$ as $k \rightarrow \infty$, we obtain that
\begin{eqnarray}\label{c1}
\bar{x} \geq 0,~~\bar{x} \neq 0.
\end{eqnarray}
If $i \in \{i \in [n] : [{\cal A}(x^k)^{m-1}+q]_i \leq 0 \}$, then
$$
[-({\cal A}(x^k)^{m-1}+q)_i]_+=-[{\cal A}(x^k)^{m-1}+q]_i.
$$
Since $\lim_{k \rightarrow \infty}\frac{q_i}{\| x^k \|} = 0$ for all $i \in [n]$, we have
\begin{eqnarray*}
0 &=& \lim_{k \rightarrow \infty}\frac{[-({\cal A}(x^k)^{m-1}+q)_i]_+}{\| x^k \|} = \lim_{k \rightarrow \infty}\frac{-[{\cal A}(x^k)^{m-1}]_i-q_i}{\| x^k \|^{m-1}}\\
&=&\lim_{k \rightarrow \infty}\frac{-[{\cal A}(x^k)^{m-1}]_i}{\| x^k \|^{m-1}} = -[{\cal A}\bar{x}^{m-1}]_i;
\end{eqnarray*}
and if $i \in \{i \in [n] : [{\cal A}(x^k)^{m-1}+q]_i \geq 0 \}$, then
\begin{eqnarray*}
0 \leq \lim_{k \rightarrow \infty}\frac{[{\cal A}(x^k)^{m-1}+q]_i}{\| x^k \|^{m-1}} = \lim_{k \rightarrow \infty}\frac{[{\cal A}(x^k)^{m-1}]_i}{\| x^k \|^{m-1}}=
[{\cal A}\bar{x}^{m-1}]_i.
\end{eqnarray*}
Combining these two situations together, we have
\begin{eqnarray}\label{c2}
[{\cal A}\bar{x}^{m-1}]_i \geq 0,~~\forall i \in [n].
\end{eqnarray}
In addition, by using (\ref{c0-1}), we have
$$
\bar{x_i}[{\cal A}\bar{x}^{m-1}]_i = \lim_{k \rightarrow \infty}\frac{x_i^k}{\| x^k \|}\frac{[{\cal A}(x^k)^{m-1}]_i}{\| x^k \|^{m-1}} = \lim_{k \rightarrow \infty}\frac{x_i^k}{\| x^k \|}\frac{[{\cal A}(x^k)^{m-1}+q]_i}{\| x^k \|^{m-1}}
\leq 0.
$$
By using (\ref{c1}) and (\ref{c2}), we have $\bar{x_i}[{\cal A}\bar{x}^{m-1}]_i \geq 0$. This, together with the above inequality, implies
\begin{eqnarray}\label{c3}
\bar{x_i}[{\cal A}\bar{x}^{m-1}]_i = 0,~~\forall i \in [n].
\end{eqnarray}
Furthermore, by combining (\ref{c1}) with (\ref{c2}) and (\ref{c3}), we obtain that $\bar{x}$ is a nonzero solution of TCP$(0, {\cal A})$. However, since every $P$ tensor is a strictly semi-positive tensor; while if ${\cal A}$ is strictly semi-positive, then TCP$(0,{\cal A})$ has a unique solution $0$ (see \cite[Theorem 3.2]{sq15}). This derives a contradiction. So the result {\bf R1} holds.

Now, suppose that the solution set of TCP$(q, {\cal A})$ is unbounded. Then there exists an unbounded solution sequence $\{ x^k \}$ of TCP$(q,{\cal A})$ such that $\| x^k \| \rightarrow \infty$ as $k \rightarrow \infty$, and for all $i \in [n], k \geq 0$, it follows that
\begin{eqnarray}\label{E-add-0}
x^k\geq 0,\quad {\cal A}(x^k)^{m-1}+q\geq 0,\quad (x^k)^T[{\cal A}(x^k)^{m-1}+q]=0.
\end{eqnarray}
Obviously, $[{\cal A}(x^k)^{m-1}+q]_i \geq 0$ implies that
$$
\frac{[-({\cal A}(x^k)^{m-1}+q)_i]_+}{\| x^k \|}\rightarrow 0\quad \mbox{\rm as}\; k\rightarrow \infty.
$$
Thus, the solution sequence $\{ x^k \}$ satisfies (\ref{c0}); and furthermore, by using the result {\bf R1}, there exist an index $i_0$ and a positive integer $k^* > 0$ such that $x_{i_0}^{k^*}[{\cal A}(x^{k^*})^{m-1}+q]_{i_0} > 0$, which is contrary to that $x_{i}^{k}[{\cal A}(x^k)^{m-1}+q]_i = 0$ for all $i \in [n]$ and $k \geq 0$. So the solution set of TCP$(q, {\cal A})$ is bounded.

{\bf Part 2}. We now show that the solution set of TCP$(q, {\cal A})$ is closed. Suppose that $\{x^k\}$ is a solution sequence of TCP$(q,{\cal A})$ and
\begin{eqnarray}\label{E-add-1}
\lim_{k\rightarrow \infty}x_k=\bar{x},
\end{eqnarray}
we need to show that $\bar{x}$ solves TCP$(q, {\cal A})$.

Since ${\cal A}x^{m-1}+q$ is continuous, by (\ref{E-add-1}) we have
\begin{eqnarray}\label{E-add-2}
\lim_{k \rightarrow \infty}[({\cal A}(x_k)^{m-1}+q]={\cal A}\left(\lim_{k \rightarrow \infty}x_k\right)^{m-1}+q={\cal A}\bar{x}^{m-1}+q.
\end{eqnarray}
Since $\{x^k\}$ is a solution sequence of TCP$(q, {\cal A})$, we have that (\ref{E-add-0}) holds. Furthermore, by using (\ref{E-add-0}), (\ref{E-add-1}) and (\ref{E-add-2}), we can obtain that
$$
\bar{x}\geq 0, \quad {\cal A}\bar{x}^{m-1}+q\geq 0,\quad \bar{x}^T({\cal A}\bar{x}^{m-1}+q)=0.
$$
So $\bar{x}$ is a solution of TCP$(q,{\cal A})$. Therefore, we obtain that the solution set is closed.

Combining {\bf Part 1} with {\bf Part 2}, we obtain that the solution set of TCP$(q,{\cal A})$ is compact. This completes the proof.
\ep

\section{Strong $P$ Tensor and Related Properties}
\label{sec:3}

In this section, we consider the question: for which kind of tensor, TCP$(q,{\cal A})$ has the GUS-property. For this purpose, we introduce a new class of tensors, called the {\it strong $P$ tensor}, which is defined as follows.
\begin{Definition}\label{def-sP}
Let ${\cal A}=(a_{i_1\cdots i_m})\in T_{m,n}$. We say that ${\cal A}$ is a strong $P$ tensor iff $F(x)={\cal A}x^{m-1}+q$ is a $P$-function.
\end{Definition}

It is well-known that a matrix $A$ is a $P$-matrix iff $F(x)=Ax+q$ is a $P$-function. Thus, the strong $P$ tensor is a generalization of the $P$-matrix from matrix to tensor. From the definitions of the $P$ tensor and the strong $P$ tensor, it is easy to see that every strong $P$ tensor must be a $P$ tensor. The strong $P$ tensor is defined with the help of the $P$-function, so an advantage of this way is that the related results and methods associated with the $P$-function can be applied to study this class of tensors.

The following result comes from \cite[Theorem 2.3]{m74a}.
\begin{Lemma}\label{lemma-PF}
Let $F:R_+^n\rightarrow R^n$ be a $P$-function, then the corresponding CP$(F)$ has no more than one solution.
\end{Lemma}

With the help of Lemma \ref{lemma-PF}, we show the following result.
\begin{Theorem}\label{thm-sPu}
Suppose that ${\cal A}\in T_{m,n}$ is a strong $P$ tensor, then TCP$(q,{\cal A})$ has the GUS-property.
\end{Theorem}

{\bf Proof.}
Since ${\cal A}$ is a strong $P$ tensor, it follows that ${\cal A}$ is a $P$ tensor. Furthermore, it follows from \cite[Corollary 3.3]{sq14b} that TCP$(q,{\cal A})$ has a solution for every $q\in \mathbb{R}^n$. Also since ${\cal A}$ is a strong $P$ tensor, it follows that ${\cal A}x^{m-1}+q$ is a $P$-function; and hence, from Lemma \ref{lemma-PF} it follows that TCP$(q,{\cal A})$ has no more than one solution. Therefore, TCP$(q,{\cal A})$ has a unique solution for every $q\in \mathbb{R}^n$, i.e., TCP$(q,{\cal A})$ has the GUS-property.
\ep

\begin{Corollary}\label{cor-sPu}
Given ${\cal A}\in T_{m,n}$ and $q\in \mathbb{R}^n$. Suppose that $F(x)={\cal A}x^{m-1}+q$ is a $P$-function, then the corresponding CP$(F)$ has the GUS-property; and $m$ must be an even number.
\end{Corollary}

In the theory of nonlinear complementarity problems, when the involved function $F$ is a $P$-function, one can only obtain that the corresponding CP$(F)$ has no more than one solution (see Lemma \ref{lemma-PF}); when the involved function $F$ is a uniform $P$-function, one can obtain that the corresponding CP$(F)$ has the GUS-property (see \cite[Corollary 3.2]{m74}). From the first result of Corollary \ref{cor-sPu}, we see that CP$(F)$ has the GUS-property when $F(x)={\cal A}x^{m-1}+q$ is a $P$-function. In addition, from the second result of Corollary \ref{cor-sPu}, we obtain that a class of functions (i.e., $F(x)={\cal A}x^{m-1}+q$ with $m$ being odd) can not be in the class of the $P$-functions.

In the following, we investigate the relationship between the $P$ tensor and the strong $P$ tensor. Recall that every strong $P$ tensor is a $P$ tensor. The following example demonstrates that the inverse does not hold.
\begin{Example}\label{exam-3}
Let ${\cal A}=(a_{i_1 i_2 i_3 i_4})\in T_{4,2}$, where $a_{1111}=1,a_{1222}=-1,a_{1122}=1,a_{2222}=1,a_{2111}=-1,a_{2211}=1$ and all other $a_{i_{1}i_{2}i_{3}i_{4}}=0$. Obviously,
$$
{\cal A}x^3=\left(\begin{array}{c}
x_1^3-x_2^3+x_1x_2^2 \\ x_2^3-x_1^3+x_2x_1^2
\end{array}\right),
$$
then
$$
x_1({\cal A}x^3)_1=x_1^4-x_1x_2^3+x_1^2x_2^2\quad\mbox{\rm and}\quad
x_2({\cal A}x^3)_2=x_2^4-x_2x_1^3+x_1^2x_2^2.
$$
We consider the following several cases:
\begin{itemize}
\item when $x_1=x_2\neq 0$, we have $x_1({\cal A}x^3)_1=x_1^4-x_1^4+x_1^4=x_1^4>0$;
\item when only one $x_i=0$ for any $i\in \{1,2\}$, then $x_j({\cal A}x^3)_j=x_j^4>0$ for $j\neq i$;
\item when $x_1>x_2>0$, we have $x_1({\cal A}x^3)_1=x_1(x_1^3-x_2^3)+x_1^2x_2^2>0$;
\item when $x_2>x_1>0$, we have $x_2({\cal A}x^3)_2=x_2(x_2^3-x_1^3)+x_1^2x_2^2>0$;
\item when $0>x_1>x_2$, we have $x_2({\cal A}x^3)_2=x_2(x_2^3-x_1^3)+x_1^2x_2^2>0$;
\item when $0>x_2>x_1$, we have $x_1({\cal A}x^3)_1=x_1(x_1^3-x_2^3)+x_1^2x_2^2>0$;
\item when $x_1>0>x_2$, we have $-x_2x_1^3>0$, and so $x_2({\cal A}x^3)_2>0$;
\item when $x_2>0>x_1$, we have $-x_1x_2^3>0$, and so $x_1({\cal A}x^3)_1>0$.
\end{itemize}
Thus, for any $x \in \mathbb{R}^2\setminus\{0\}$, there exists an index $i\in \{1,2\}$ such that $x_i({\cal A}x^3)_i>0$. So ${\cal A}$ is a $P$ tensor. However, ${\cal A}$ is not a strong $P$ tensor. In fact, if we take $x=(2.1,-1.9)^T$ and $y=(2,-2)^T$, then we have
$$
(x_1-y_1)(({\cal A}x^3)_1-({\cal A}y^3)_1)=-0.0299<0
$$
and
$$
(x_2-y_2)(({\cal A}x^3)_2-({\cal A}y^3)_2)=-0.0499<0.
$$
Thus, by Definition \ref{def-sP} we obtain that ${\cal A}$ is not a strong $P$ tensor.
\end{Example}

Example \ref{exam-3} demonstrates that the set consists of all strong $P$ tensors is a proper subset of the set consists of all $P$ tensors.

Many properties of the $P$ tensor have been obtained in the literature. Since every strong $P$ tensor is a $P$ tensor, we may easily obtain the following properties of strong $P$ tensor:
\begin{Proposition}\label{prop-sP}
If ${\cal A} \in T_{m,n}$ is a strong $P$ tensor, then
\begin{itemize}
\item[(i)] ${\cal A}$ must be strictly semi-positive;
\item[(ii)] ${\cal A}$ must be an $R$ tensor;
\item[(iii)] all of its H-eigenvalues and Z-eigenvalues are positive;
\item[(iv)] all the diagonal entries of ${\cal A}$ are positive;
\item[(v)]every principal sub-tensor of ${\cal A}$ is still a strong $P$ tensor.
\end{itemize}
\end{Proposition}

{\bf Proof.}
Since a strong $P$ tensor is a $P$ tensor, the first four results can be easily obtained from \cite{sq14a,sq14b}. Now we prove the result (v). Let an arbitrary principal sub-tensor ${\cal A}_r^J\in T_{m,r}$ of the strong $P$ tensor ${\cal A}\in T_{m,n}$ be given. We choose any $x=(x_{j_1}, x_{j_2}, \cdots, x_{j_r})\in \mathbb{R}^r\setminus\{0\}$ and $y=(y_{j_1}, y_{j_2}, \cdots, y_{j_r})\in \mathbb{R}^r\setminus\{0\}$ with $x\neq y$. Then let $\bar{x}=(\bar{x}_1, \bar{x}_2, \cdots, \bar{x}_n) \in \mathbb{R}^n$ where $\bar{x}_i=x_{j_i}$ for $i \in J$ and $\bar{x}_i=0$ for $i\notin J$. In a similar way, let $\bar{y}=(\bar{y}_1, \bar{y}_2, \cdots, \bar{y}_n) \in \mathbb{R}^n$ where $\bar{y}_i=y_{j_i}$ for $i \in J$ and $\bar{y}_i=0$ for $i \notin J$. Since ${\cal A}$ is a strong $P$ tensor, then there exists an index $k \in [n]$ such that
\begin{eqnarray*}
0&<& \max_{k\in [n]}\;(\bar{x}_k-\bar{y}_k)(({\cal A}\bar{x}^{m-1})_k-({\cal A}\bar{y}^{m-1})_k)\\
&=&\max_{k\in J}\;(x_k-y_k)(({\cal A}_r^Jx^{m-1})_k-({\cal A}_r^Jy^{m-1})_k).
\end{eqnarray*}
Thus, ${\cal A}_r^J$ is a strong $P$ tensor.
\ep

\section{Conclusions}
\label{sec:4}

By constructing two counter-examples, we proved that ${\cal A}$ is a $P$ tensor does not imply that TCP$(q,{\cal A})$ has the GUS-property; and that TCP$(q,{\cal A})$ has the GUS-property does not imply that ${\cal A}$ is a $P$ tensor. These gave a negative answer to Question 6.3 proposed in \cite{sq14a}. We also showed that the solution set of TCP$(q,{\cal A})$ is nonempty and compact when ${\cal A}$ is a $P$ tensor.

In order to investigate that for which kind of tensor, the tensor complementarity problem has the GUS-property, we introduced the concept of the strong $P$ tensor, and showed that TCP$(q,{\cal A})$ has the GUS-property when ${\cal A}$ is a strong $P$ tensor. We also proved that every strong $P$ tensor is a $P$ tensor; and hence, many known results associated with the $P$ tensor were generalized to the case of the strong $P$ tensor.

Note that the strong $P$ tensor is defined by using the $P$-function, we believe that more properties related to the strong $P$ tensor can be further studied with the help of known methods and results for the $P$-function.




\end{document}